\documentstyle[12pt]{article}
\newcommand{\const}{\mathop{\rm const}\limits}

\begin{document}

\begin{center}

{\bf  MIXED LEBESGUE SPACE NORM STRICHARTZ TYPE  \\

\vspace{3mm}

ESTIMATION  FOR  SOLUTION OF INHOMOGENEOUS \\

\vspace{3mm}

 PARABOLIC EQUATION, \\

\vspace{3mm}

with constants  evaluation. } \par

\vspace{4mm}

 $ {\bf E.Ostrovsky^a, \ \ L.Sirota^b } $ \\

\vspace{4mm}

$ ^a $ Corresponding Author. Department of Mathematics and computer science, Bar-Ilan University, 84105, Ramat Gan, Israel.\\
\end{center}
E - mail: \ galo@list.ru \  eugostrovsky@list.ru\\
\begin{center}
$ ^b $  Department of Mathematics and computer science. Bar-Ilan University,
84105, Ramat Gan, Israel.\\

E - mail: \ sirota3@bezeqint.net\\

\vspace{5mm}
                    {\bf Abstract.}\\

 \end{center}

 \vspace{4mm}

 We give an optimal in mixed (anisotropic)  Strichartz type Lebesgue space-time norm  estimates for the solution of linear
parabolic inhomogeneous initial problem, with are exact or exact up to multiplicative constant coefficient evaluation. \par

\vspace{4mm}

 \par

\vspace{4mm}

  \par

\vspace{4mm}

\section{Notations. Statement of problem.}

\vspace{3mm}

{\bf Statement of problem.} \par

\vspace{3mm}

We  consider in this article the classical initial value problem  for the function $ u = u(x,t), \ t \in T = (0,I), \ 0 < I \le \infty,  $
of a multivariate linear Parabolic PDE equations in whole Euclidean space $ x \in X = R^d  $ with  density of external force   $ F(x,t) $ and
initial condition $  f(x) $ of a form

$$
\partial_t u = 0.5 \ \Delta u + F(t,x), \ u(x, 0+) = f(x). \eqno(1.1)
$$

 We will update after the statement of this problem: for given $  p \in (1, \infty) $

 $$
 \lim_{ t \to 0+} | u(\cdot,t) - f(\cdot)|L_p(R^d) = 0.
 $$

  Let us introduce the following notations.\\

$$
x = \vec{x} = \{x_1, x_2, \ldots, x_d \} \in R^d \ \Rightarrow  |x| = \sqrt{(x,x)} = \sqrt{\sum_{i=1}^d x_i^2},
$$

$$
w_t(x) = (2 \pi t)^{-d/2} \ \exp \left( - \frac{|x|^2}{2t}  \right), \ t > 0.
$$
 This function is the fundamental solution of  classical heat equation ("heat potential").\par

Further, we denote as usually the convolutions

$$
F*G(x) := \int_X F(x-y) \ G(y) \ dy, \ x,y \in R^d;  \
$$
("space" convolution); and at the same notation will be used  for a "time" convolution:

$$
f*g(t) := \int_0^t f(t-s) \ g(s) \ ds, \ s,t \ge 0.
$$
 Authors hope that this definitions not lead to the confusion. \par

  The unique "regular" solution of the problem (1.1) may be written under known natural conditions
 imposed on the data $ f, \ F $  (measurability, boundedness, belonging to and or other Banach space etc.)
 in explicit view  as follows; $  u(x,t) = u_0(x,t) + u_1(x,t), $ where

  $$
  u_0(x,t) = w_t* [ f](x) = \int_X w_t(x-y) \ f(y) \ dy,  \eqno(1.2)
  $$

$$
u_1(x,t) = [w_t ** F](x,t) = \int_0^t ds \int_X w_{t-s}(x-y) \ F(y,s) \ dy. \eqno(1.3)
$$

 \vspace{3mm}

 Further, denote  as ordinary

$$
|f|_r = |f|_{r,X} = \left[ \int_X |f(x)|^r \ dx  \right]^{1/r}, \ r = \const \ge 1;
$$

$$
|g|_{q} = |g|_{q,T} =  \left[\int_T |g(t)|^q \ dt  \right]^{1/q}, \ q = \const \ge 1;
$$

The so - called {\it mixed} $ L_{p,q} = L_{p,X; q,T}  $ norm  of a function of "two" variables
 $ F = F(x,t) $ is  defined by equality:

 $$
  |F|_{p,X; q,T} := \left\{\int_0^{\infty} \left[ \int_X |F(x,t)|^p \ dx \right]^{q/p} \ dt \right\}^{1/q}, \ 1 \le p,q < \infty,
  \eqno(1.4)
 $$

 and analogously

  $$
 |F|_{q,T; p,X}:= \left\{\int_X \left[ \int_0^{T} |F(x,t)|^q \ dt \right]^{p/q} \ dt \right\}^{1/p}, \ 1 \le p,q < \infty,
  \eqno(1.4a)
 $$

  Not to be confused with the standard Lorentz norm! \par
 The space consisting on all the (common measurable) function $  F(x,t) $  with finite mixed norm $ | F(\cdot, \cdot) |_{p,q}  $ is
 said to be {\it  anisotropic, } or equally {\it Bochner} space  and denoted similar  $ L_{p,q} = L_{p,X; q,T}.  $ \\

  These spaces appear in an article of   Benedek A. and Panzone R. \cite{Benedek2} and was completely
investigated  in the classical monograph  of  Besov O.V., Ilin V.P., \newline Nikolskii S.M. \cite{Besov1}.\par
 It is known, see  \cite{Adams1}, \cite{Besov1}, p. 7, 25, that if $ p_1 \le p_2  $

$$
|F|_{p_1, p_2} \le |F|_{p_2, p_1},
$$
"permutation inequality". \par

Note that in general case $ |F|_{p,q} \ne |F|_{q,p}, $ but $ |F|p,p = |F|p. $ \par

 Observe also that if $ F(x, t) = G_1(x) \cdot G_2(t), $  (condition of factorization), then
$ |F|_{p,q} = |G_1|p \cdot |G_2|q, $ (formula of factorization).\par

 These spaces arise in the Theory of Approximation, Functional Analysis, theory
of Partial Differential Equations, theory of Random Processes etc.   Consider for instance the linear
integral operator $U $ acting on the functions defined on the measurable space $ (Y, B, \nu) $ onto another
measurable space $ (X, A, \mu) $ with common  measurable kernel $ K(x,y):  $

$$
U[f](x) = \int_Y K(x,y) \ f(y) \ \mu(dy).
$$
then we conclude by means  of H\"older's inequality

$$
|U[f]|_{p,X} \le |K|_{p,q'} \cdot |f|_{q,Y}; \ q' \stackrel{def}{=} q/(q-1), \ 1 < q \le \infty,
$$
or equally

$$
|U|_{q \to p} \stackrel{def}{=}   \sup_{f \ne 0, \ f \in L_q(Y)} \left[ \frac{|U[f]|_{p,X}}{|f|_{q,Y}} \right] \le |K|_{p,q'}.
$$

\vspace{4mm}

{\bf  Our purpose is to estimate both the functions  $ u_0 $ and $ u_1  $ in the mixed norms $ |u_0|_{r,T;q,X}, \  |u_1|_{p,X;q,T}  $
via correspondingly also mixed norms of the data  $ | f|_{r,X}   $ and $ |F|_{r,X; k,T}.  $ } \par

\vspace{4mm}

 {\it We simplify known proofs and estimates, write down the asymptotically optimal  estimates of the constants and prove its
 exactness.    } \par

\vspace{4mm}

Previous works:  \cite{Benedek1}, \cite{Krylov1} - \cite{Krylov3}, \cite{Kyeong1}, \cite{Ladyzhenskaya2} - \cite{Ladyzhenskaya1},
  \cite{Zhai1} etc. The applications of these estimates in the theory of non - linear evolutionary  PDE, for example, equations
 of Navier - Stokes see in  \cite{Cui1} -  \cite{Zhang1}.  \par

\vspace{4mm}

 The paper is organized as follows. It the second section we obtain the mixed norm estimates for the first component of solution $ u_0= u_0(x,t), $
 in the third  we obtain ones for the second component   $ u_1 = u_1(x,t). $ \par
  In the fourth section we will prove the non-refinements of our conditions, where we prove simultaneously the exactness of conditions in the
 Young inequality convolution. The next section is devoted to common, or equally united estimates.
 We consider in the sixth section some generalizations of obtained estimates into the so-called Grand Lebesgue Spaces. \par
 The $ 7^{th} $ section contains the weight estimates for the classical inhomogeneous parabolic initial value problem. The penultimate
section described the mixed norm estimates for solution of non-local evolution initial value problem with fractional Laplace's operator. \par

 By tradition, the last section contains some concluding remarks. \par

\section{Influence of initial condition.}

\vspace{4mm}

We intend to obtain the estimation of a form

$$
 |u_0|_{r_0,T; q_0,X} =  |w_t*f|_{r_0,T; q_0,X} \le K_0(d;p_0,q_0, r_0) \ |f|_{p_0}. \eqno(2.1)
$$

 {\it Notations and restrictions: } \\

$$
p_0,q_0,r_0 \in (1,\infty), \hspace{5mm}  d \ r_0 > 2; \  z = |x-y| > 0; \ q_0 \in (r_0 d/(r_0 d - 2), \infty);
$$

$$
B(z) = B_r(z) := |w_t(z)|_{r,T} = (2 \pi)^{-d/2} \
\sqrt[r]{\int_{R_+} t^{-dr/2} \ e^{ -r \ z^2/2t } \ dt} =
$$

$$
 \pi^{d/2} \ z^{2/r-d} \ r^{1/r - d/2} \ \Gamma^{1/r}(dr/2 -1) := D(d,r) \ |x-y|^{2/r - d},
$$
where

$$
D(d,r) = \pi^{d/2} \ r^{1/r - d/2} \ \Gamma^{1/r}(dr/2 -1).
$$

 Note that as $  r \to d/2 + 0  $

$$
D(d,r) \sim (2 \pi)^{d/2} \ (rd - 2)^{-d/2}.
$$

 Further, put

$$
V(p,r) := \frac{(d - 2/r)^{-1}}{ \left[(p-1)(dr/2 - p) \right]^{\kappa} };
$$

$$
\kappa: = 1 - \frac{2}{d r} \in (0,1); \ 1 + \frac{1}{q} = \frac{1}{p} + \kappa \ \Leftrightarrow \ \frac{1}{q} = \frac{1}{p} - \frac{2}{rd};
\eqno(2.2)
$$

 Let us introduce the following {\it family } of domains $ G_0(p) = G_0(d;p), \ p > 1 $ on the plane $ (q,r): $

$$
G_0(p) = G_0(d; p) = \{ (q,r): q > 1, \ r > 2/d,  \ 1/q + 2/dr = 1/p \}.
$$

\vspace{4mm}

{\bf Theorem 2.0a.}  {\it Let $ f \in L_{p_0}(R^d) $ for some $ p_0 > 1. $ There holds  under described restrictions:

$$
(q_0,r_0) \in G_0(d; p_0)
$$
 the following estimate:}

$$
|u_0|_{r_0,T; q_0,X} \le  C_0(d; p_0) \ D(d,r_0) \ V(p_0,r_0) \ |f|_{p_0}. \eqno(2.3)
$$

{\bf Remark 2.1.} We can adopt as the capacity of the value $   C_0(d; p_0) $  its minimal value, indeed:

 $$
 C_0(d; p) :=  \sup_{ (q,r) \in G_0(p) } \sup_{f: |f|_p = 1} \
 \left[ \frac{|u_0|_{r,T; q,X}}{D(d,r) \ V(p,r) \ |f|_p }  \right] < \infty, \ 1 < p < \infty.
 $$

  \vspace{4mm}

 {\bf Proof.} We can and will  suppose without loss of generality that $  f(x) \ge 0. $   We deduce using
 Marcinkiewicz integral triangle  inequality  writing locally $  (p,q,r)  $   instead $ (p_0,q_0,r_0):  $
 $$
 |u_0|_{r,T} \le \int_X |w_t(x-y)|_{r,T} \ f(y) \ dy = D(d,r) \int_X \frac{f(y) \ dy}{|x-y|^{d - 2/r}}.
$$

  Notice that the integral in the right - hand side is the well -  known  fractional Riesz potential:

$$
 \int_X \frac{f(y) \ dy}{|x-y|^{d - 2/r}} = \int_{R^d} \frac{f(y) \ dy}{|x-y|^{d - 2/r}} = I_{2/r}[f](x) = I^{(d)}_{2/r}[f](x).
$$

   The using for us exact
up to multiplicative constants  Lebesgue - Riesz estimates for this operator are obtained   in the article
\cite{Ostrovsky108}:

$$
|I^{(d)}_{2/r}[f](\cdot)|_{q.X} \le C_0(p,r) \ V(p,r) |f|_p, \ p,q,r \in (1,\infty), \ \frac{1}{q} = \frac{1}{p} - \frac{2}{d r}. \eqno(2.4)
$$

 This completes the proof of theorem 2.0. \par

\vspace{3mm}

{\bf Remark 2.2.}   The condition (2.2) coincides with  ones in proposition 1.2 in the articles \cite{Miao1}, \cite{Miao2}, \cite{Zhai1},
where is consider also the case of {\it fractional } Laplace operator $  (-\Delta)^{\alpha} $ instead the classical operator $  (- \Delta), $
but without constants estimation. \par
 The method used in \cite{Zhai1} is different on our way; for instance, we do not use the theory of interpolation of operators. \par

\vspace{3mm}

{\bf Theorem 2.0b.}  {\it Let $ f \in L_{p_0}(R^d) $ for some $ p_0 > 1. $ There holds  the following estimate:}

$$
|u_0|_{q_0,X} \le   K_W(m) \  K_B(m,p_0) \ |f|_{p_0,X} \ t^{-d/2(1/p_0-1/q_0)}, \eqno(2.5)
$$
{\it where }

$$
q_0 \ge p_0, \hspace{6mm}  1 + 1/q_0 = 1/m + 1/p_0, \ q, p_0,m > 1. \eqno(2.6)
$$

{\it Moreover, the relation (2.6) is necessary for the inequality of a form (2.5). }\par
 {\it Besides, the asymptotical equality in (2.5) as $ t \to \infty $
  is attained iff  $ f_0  $  is  density of Gaussian centered non - trivial  distribution:  } \par

$$
f_0(x) = ( 2 \pi)^{-d/2} \ \sigma^{-d} \ \exp \left(-|x|^2/(2 \sigma^2) \right), \sigma = \const > 0.
$$

\vspace{3mm}

{\bf Proof.} We have: $ u_0(x,t) = w_t*[f](x).$ We apply the Young's - Beckner's inequality:

$$
|u_0|_{q_0,X} \le K_B(d; m,p_0) \ |w_t|_{m,X} \ |f|_{p_0,X} = K_B(d; m,p_0) \ K_w(m) \ |f|_{p_0} \ t^{-d/2(1/p_0 - 1/q_0)}.
$$

 The necessity of equality (2.6) may be easily proved by means of scaling method. The last proposition about
asymptotical equality  as $ t \to \infty $  is in fact proved  by  W.Beckner in \cite{Beckner1}. \par

\vspace{3mm}

 Note that the inequality (2.5) is well-known, see for example, \cite{Kato1}; we write only the constants estimates. \par

\vspace{4mm}

\section{Influence of right - hand side.}

\vspace{4mm}

 Some new notations and restrictions (conditions).  Let four numbers $ (p,q,r,k) \in (1,\infty)  $ be a given. Denote

 $$
 Q = Q(p,r): \ \frac{1}{p} + 1 = \frac{1}{Q} + \frac{1}{r};  \hspace{4mm} \theta = \theta(p,r) = \frac{d}{2} \left( 1 - \frac{1}{Q} \right),
 $$

$$
k_- = 1, \ k_+ = \frac{d}{d - \theta}, \hspace{5mm} q_- = \frac{q}{\lambda}   = \frac{2}{1/r - 1/p}, \ q_+ = \infty,
$$

$$
\beta = \frac{\theta}{d} =   \frac{1}{2} \left(  \frac{1}{r} - \frac{1}{p}   \right).\eqno(3.0)
$$

 Further, we denote  the  $  L_Q $ norms  of heat potential:
$$
K_w = K_w(d,Q):=  |w_1|_{Q} =  |w_1|_{Q,X} = (2 \pi )^{d(1-Q)/(2 Q)}  \ Q^{- d/2Q };  \ Q \ge 1.
$$
 As a consequence:
$$
|w_t|_{Q} = (2 \pi )^{d(1-Q)/(2 Q)} \ t^{0.5 d(1/Q - 1)} \ Q^{- d/2Q } = K_w(d,Q) \ t^{0.5 d(1/Q - 1)}, \hspace{4mm} \ t > 0.
$$

\vspace{3mm}

 We introduce also the following {\it family } of domains $ G_1(k,r) = G_1(d;k,r), \ k,r > 1 $ on the plane $ (p,q): $

$$
G_1(k,r) = G_1(d;k,r) = \{ (p,q):   \ p > r, \ q > 1, \ 1/q + d/2p = 1/k + d/(2r) -1 \}. \eqno(3.1)
$$

\vspace{3mm}

{\bf Theorem 3.1a.} {\it  Let again the four numbers $ (p_1,q_1,r_1,k_1) \in (1,\infty)  $ be a given. }
{\it Suppose} $  F (\cdot, \cdot) \in L(r_1,X; k_1,T) $  {\it for some $ k_1,r_1 > 1. $
If  $  (q_1,r_1) \in G_1(d; k_1,r_1), $ } {\it and } $ k_1 \in (1, \ k_+), $ {\it then}

$$
|u_1(\cdot, \cdot)|_{p_1,X; q_1,T} \le  K_1(p_1,q_1,r_1,k_1) \  |F (\cdot, \cdot)|_{r_1,X; k_1,T}, \eqno(3.2)
$$
{\it  where}

$$
K_1(p_1,q_1,r_1,k_1) = \frac{C_1(d; r_1,k_1)}{ \left[ k_1-k_- \right]^{\beta}}. \eqno(3.3)
$$

 {\it  Conversely, the equality  $ (q_1,r_1) \in G_1(d; k_1,r_1) $ is necessary for the estimation of the form (3.2). }\par
{\it  Moreover, the estimation (3.2) is  exact up to multiplicative constant. } \par

\vspace{3mm}

{\bf Remark 3.1.} We can adopt as before as the capacity of the value $  C_1(d; r,k) $  its minimal value, namely:

$$
C_1(d; r,k) := \sup_{ (p,q) \in G_1(k,r) }
 \left[ \frac{\left[ k-k_- \right]^{\beta} \ |u_1(\cdot, \cdot)|_{p,X; q,T}}{|F (\cdot, \cdot)|_{r,X; k,T}} \right] < \infty, \ r,k > 1.
$$

 \vspace{3mm}

 {\bf Remark 3.2.}  Note that the norms in the proposition (3.2) follows in reverse order as in the assertion of
 theorem 2.0. \par

 \vspace{3mm}

 {\bf Proof.}  We get using again Marcinkiewicz (triangle) inequality writing again $ (p,q,r,k)  $ instead $ (p_1,q_1,r_1,k_1):  $

 $$
 |u_1(\cdot, \cdot)|_{p,X} = \left|\int_0^t ds \int _X w_{t-s}(x-y) \  F(y,s) \ dy  \right|_{p,X} \le
 $$

 $$
 \int_0^t ds \int _X  \left| w_{t-s}(x-y) \ F(y,s) \right|_{p,X} \ dy = \int_0^t \left|w_{t-s}(\cdot)*F(\cdot,s) \right|_{p,X} \ ds.
 $$
  We apply the Young's inequality for the space convolution:

 $$
 |u_1(\cdot, \cdot)|_{p,X} \le K_B(d; Q,r) \int_0^t |w_{t-s}(\cdot)|_{Q,X} \ \cdot |F(\cdot,s)|_{r,X} \ ds =
 $$

 $$
   K_B(d; Q,r) \ K_w(d,Q) \int_0^t (t-s)^{0.5 d(1/Q - 1)} \ |F(\cdot, s)|_{p,X} ds, \eqno(3.4)
 $$
  where $ K_B(d; p,q) < 1 $ is the famous Beckner's constant:

$$
A(p):=  \left[ \frac{p^{1/p}}{p'^{1/p'}} \right]^{1/2}, \ p > 1, \ p' := p/(p-1);
$$

$$
K_B = K_B(p,q) = K_B(d; p,q) :=  \left( A(p) \ A(q) \ A(r) \right)^d, \ 1 + 1/r = 1/p + 1/q.
$$

 W.Beckner in \cite{Beckner1} proved that

$$
|f * g|_r \le K_B(d; p,q) \ |f|_p \ |g|_q,
$$
where again

$$
1 + \frac{1}{r} = \frac{1}{p} +  \frac{1}{q}. \eqno(3.5)
$$

 See also \cite{Brascamp1}.\\

 The integral in the right - hand (3.4) is so-called Rieman's fractional integral or equally  Cesaro-Hardy average
operator. The Lebesgue - Riesz $ L_q(T)  $ norm estimates    for this operator  with exact up to multiplicative constants
 evaluating  are obtained in \cite{Ostrovsky109}, see also \cite{Ostrovsky108}.  Indeed, if

 $$
 U_{\theta}[h](\cdot) = \int_0^t (t-s)^{-\theta} \ h(s) \ ds, \ \theta = \const \in (0,1),
 $$
 then

 $$
 |  U_{\theta}[h](t) |_{q,T} \le K \cdot |h|_{k,T},
 $$
 where

 $$
 1 + \frac{1}{q} =  \frac{1}{k} + \theta.
 $$
  Here we must substitute

  $$
  \theta = \frac{d}{2} \left( 1 -   \frac{1}{Q} \right) \in (0,1).
  $$

  It remains to use this estimations. \par
 This completes the proof of theorem 3.1a.\par

\vspace{3mm}

{\bf Remark 3.3.}   The condition (3.1)  coincides with  ones
   in equality (1.8) in the article  \cite{Zhai1}, where is consider also the case
of {\it fractional } Laplace operator $  (-\Delta)^{\alpha} $ instead the classical operator $  (- \Delta), $  but
{\it without constants estimation.} \par

The method used in \cite{Zhai1} is different on our way. \par

\vspace{3mm}

 We  investigate now the inverse order of the powers $ (p,q). $

\vspace{3mm}

{\bf Theorem 3.1b.} {\it  Let the constants  }  $ (h,k,q,r) $  {\it be such that }

$$
1 < h,k,q,r  < \infty, \hspace{6mm} \frac{1}{m} =  \frac{1}{h}- \frac{2}{dr},  \hspace{6mm}  1 +  \frac{1}{q} = \frac{1}{r} + \frac{1}{k}. \eqno(3.6)
$$

 {\it Then }

$$
|u_1|_{q,T; m,X} \le D(d,r) \ K_B(1; r,k) \ C_0(h,r) \ V(h,r) \ |F|_{k,T;h,X}, \eqno(3.7)
$$

 {\it  Conversely, the equalities   (3.6) is also necessary for the estimation of the form (3.7). }\par
{\it  Moreover, the estimation (3.7) is  exact up to multiplicative constant. } \par

\vspace{3mm}

{\bf Proof } is alike to ones in theorem 2.0a. We use  again the Marcinkiewicz's and Young's-Beckner's inequalities:

$$
|u_1|_{q,T} \le \int_X dy \ \int_0^t | w_{t-s}(x-y) \ F(y,s) |_{q,T} \ ds \le
$$

$$
K_B(1;r,k) \ \int_X |w_{\cdot}(x-y)|_{r,T} \ |F(y,\cdot)|_{k,T} \ dy =
$$

$$
K_B(1;r,k) \ D(d,r) \int_X \frac{|F(y,\cdot)|_{k,T} \ dy}{|x-y|^{d-2/r}} = K_B(1;r,k) \ D(d,r) \ I_{2/r}[F(y,\cdot)|_{k,T}].
$$

 It remains to use the estimate (2.4) for Riesz's fractional potential.\par

 The last propositions of considered theorem will be proved further, in the section 8.\par

\vspace{4mm}

\section{Necessity  of our conditions.}

\vspace{4mm}

 We intend to prove in this section that our estimates are essentially non - improvable. \par

 \vspace{3mm}

  We will use the so - called {\it dilation,  } or equally {\it  scaling method, } see \cite{Stein1}, chapter 9,
  \cite{Talenti1}. Recall the definition: the operator  $  T_{\lambda}[f], \ \lambda \in (0,\infty), $ (more exactly,
  the family of operators)  of a form

 $$
  T_{\lambda}[f](x) = f(\lambda \ x), \ x \in X = R^d \eqno(4.0)
 $$
 is said to be the {\it dilation operator.}  Here $  f(\cdot) $ belongs to a certain class of functions, for instance,
 $  L_p(X) $ space or Schwartz set $ S(R^d) $ etc. Obviously, if $ f \in   S(R^d),  $ then  $ T_{\lambda}[f] \in S(R^d). $   \par

\vspace{4mm}

{\bf 1. Young's inequality. } \par

\vspace{3mm}

 We begin our considerations of this section from the famous Young's inequality for convolution. Namely,

 $$
 |f*g|_r \le K_B( p,q) \ |f|_p \ |g|_q, \ p,q,r \in (1,\infty),  \eqno(4.1)
 $$
where $ K_B( p,q) = K_B(d; p,q) < 1 $ is also the Beckner's constant.\\

\vspace{3mm}

{\bf Theorem 4.1.}  {\it Suppose that there exists a finite constant $ K = K_B(d; p,q) $ for which
the inequality (4.1) is satisfied for arbitrary pair of functions $ (f,g) $
from the Schwartz space $ S(R^d). $  Then the triple $ (p,q,r) $ satisfies the equality (3.5). } \par

\vspace{3mm}

{\bf Proof.} Assume the inequality (4.1) is satisfied  for any  pair of functions from the Schwartz  space; then for arbitrary
number  $  \lambda \in (0, \infty)  $ and some $ f,g \in S(R^d), \ f, g \ne 0 $

 $$
 |T_{\lambda}f*  T_{\lambda} g|_r \le K_B(d; p,q) \ |T_{\lambda} f|_p \ |T_{\lambda} g|_q, \ p,q,r \in (1,\infty),  \eqno(4.3)
 $$
 Note that

 $$
 |T_{\lambda} g|_q = \lambda^{-d/q}|g|_q, \hspace{5mm} |T_{\lambda} f|_p =  \lambda^{-d/p}|f|_p,
 $$

$$
|T_{\lambda}f*  T_{\lambda} g|_r = \lambda^{-d - d/r } \ | f*g|_r.
$$
 We deduce substituting into (4.3)

$$
\lambda^{-d - d/r } \le C \cdot \lambda^{-d/p - d/q}, \ C \in (0,\infty), \eqno(4.4)
$$
 where $  C  $ does not dependent on the variable $ \lambda. $ \par
 Since the number $ \lambda $ is arbitrary positive, we conclude  from (4.4) $ 1 + 1/r = 1/p + 1/q, $ Q.E.D. \\

\vspace{3mm}

 The multivariate version of Young's inequality  has a form, see \cite{Besov1}, chapter 2, formula (2.14):

$$
| f*g|_{\vec{r}} = \left| \int_{R^d} f( x - y) \ g(y) \right|_{\vec{r}}  \le |f|_{\vec{p}}  \ |g|_{\vec{q}},
$$
where

$$
\forall j = 1,2,\ldots,d \ \Rightarrow 1 + \frac{1}{r_j} = \frac{1}{p_j} + \frac{1}{q_j}
$$
and wherein the last relation is necessary for the multivariate version of Young's inequality.\par
 This proposition may  be proved by means of consideration of the so-called factorable functions

$$
f(x) = \prod_{j=1}^d f_j(x_j), \hspace{5mm} g(x) = \prod_{j=1}^d g_j(x_j).
$$

\vspace{4mm}

{\bf 2. Initial condition.} \\

\vspace{3mm}

{\bf Theorem 4.2.}  {\it Suppose that there exists a finite constant $ K^{(0)} = K^{(0)}_B(d; p,q) $ for which
the inequality (2.3) is satisfied for arbitrary  functions $ f $ from the Schwartz space $ S(R^d): $

$$
|u_0|_{r,T; q,X} \le  K^{(0)}_B(d; p,q)  \ |f|_p, \ p,q,r = \const \in R.  \eqno(4.5)
$$

 Then the triple $ (p,q,r) $ satisfies the relation  }

$$
\frac{1}{q} =  \frac{1}{p} -  \frac{2}{d r}. \eqno(4.6)
$$

\vspace{3mm}

{\bf Proof } is completely analogous to ones in the theorem 4.1.
 Suppose that the inequality (4.5) is valid  for some  functions $ f \ \ne 0 $ from the Schwartz  space; then for arbitrary
number  $  \lambda \in (0, \infty)  $

 $$
 |  w_t * T_{\lambda}f|_{r,T; q,X} \le K_B^{(0)}(d; p,q) \ |T_{\lambda} f|_p.  \eqno(4.7)
 $$
 Note that

 $$
 w_t(x/\lambda) = \lambda^d w_{t \lambda^2}(x), \ t > 0, \ x \in R^d, \eqno(4.8)
 $$
  and we find as before after simple computations using identity (4.8)

 $$
 |  w_t * T_{\lambda}f|_{r,T; q,X} =  |  w_t * f|_{r,T; q,X} \cdot   \lambda^{-d/q - 2/r}. \eqno(4.9)
 $$

 We deduce substituting into (4.7)

$$
\lambda^{-d/q - 2/r } \le C_0 \cdot \lambda^{-d/p }, \ C_0 \in (0,\infty), \eqno(4.10)
$$
where $  C_0  $ does not dependent on the variable $ \lambda. $  We conclude equating the exponents by $ \lambda: $

$$
-d/q - 2/r  = -d/p, \ \Leftrightarrow 1/q = 1/p - 2/(d r).
$$

\vspace{3mm}

{\bf 3.  Right - hand side.} \hspace{5mm} Analogously may be proved the following result. \\

\vspace{2mm}

{\bf Theorem 4.3.}   {\it Suppose  that there exists finite constant $ K^{(1)} = K^{(1)}_B(d; p,q,r,k) $ for which
the inequality (3.1) is satisfied for arbitrary functions $ F = F(x,t) $ from the Schwartz space $ S(R^{d+1}): $

$$
|u_1(\cdot, \cdot)|_{p,X; q,T} \le  K^{(1)}_B(d; p,q,r,k) \  |F (\cdot, \cdot)|_{r,X; k,T}. \eqno(4.11)
$$

 Then the tuple of the four numbers $ (p,q,r,k) $ satisfies the equality  }

$$
1 + \frac{1}{q} - \frac{1}{k} = \frac{d}{2} \left(  \frac{1}{r} -  \frac{1}{p} \right). \eqno(4.12)
$$

\vspace{3mm}

{\bf Proof.} We can use for convenient as a capacity of the test function  $  F  $ the  factorable expression:

$$
F_0(x,t)  = g(x) \ h(t) \stackrel{def}{=} [g \otimes h](x,t). \eqno(4.13)
$$

 We introduce for these functions a following (linear) operator, more exactly, the family of operators (generalized dilation):

$$
S_{\lambda}  [g \otimes h](x,t) := g(\lambda x) \cdot h(\lambda^2 t) = T_{\lambda}[g](x) \cdot T_{\lambda^2}[h](t).
$$
 If the inequality (4.11)  is true  for arbitrary non - zero functions $  g,h $ from the correspondent Schwartz space, then

$$
|W_{t}** S_{\lambda} F_0|_{p,X; q,T} \le  K^{(1)}_B(d; p,q,r,k) \  |T_{\lambda} g|_{r,X} \cdot |T_{\lambda^2} [h]_{k,T}. \eqno(4.14)
$$

 We have consequently:

$$
|T_{\lambda} g|_{r,X} \cdot |T_{\lambda^2} [h]_{k,T} = |g|_{r,X} \ |h|_{k,T} \ \lambda^{- d/r - 2/k},
$$

$$
|W_{t}** S_{\lambda} F_0|_{p,X; q,T} = |W_{t}** F_0|_{p,X; q,T} \ \lambda^{-2 - d/p - 2/q }.
$$
 We conclude as before  after substituting into (4.14):

$$
-2 - d/p - 2/q = - d/r - 2/k.
$$
 This completes the proof of theorem 4.3. \\

\vspace{4mm}

\section{United estimations.}

\vspace{4mm}

 Recall that  the solution $ u = u(x,t) $ of the source equation   (1.1) with correspondent initial condition may be
represented as a sum  $ u(x,t) = u_0(x,t) + u_1(x,t). $  We deduce as a slight consequence synthesizing the assertions of
theorems 2.0a and 3.1b: \par

\vspace{3mm}

{\bf  Proposition 5.1.} {\it  Let the tuple   }  $ (q_1, m_1, k_1, h_1)  $ {\it  satisfies simultaneously the conditions (2.2) and (3.6). } \par
{\it Let also } $ f \in L_{p_0}(R^d), \ F \in L_{ k_1,T; h_1,X}((0,T) \otimes R^d).  $ {\it Then }

$$
|u|_{q_1,T; m_1.X} \le C_0(d,p_0) \ D(d,q_1) \ V(p_0,q_1) \ |f|_{p_0} +
$$

$$
C_0(h_1,r_1) \ D(d,r_1) \ V(h_1, r_1) \ K_B(1;r_1,k_1) \ |F|_{k_1,T;h_1,X}. \eqno(5.1)
$$

 Let us consider more general case.  Recall previously that the direct sum

  $$
M_{r_0, q_0, q_1, m_1} = M_{r_0, q_0, q_1, m_1}( (0,T) \otimes R^d)  =
  $$

  $$
 L_{r_0, q_0}( (0,T) \otimes R^d)  \oplus  L_{q_1, m_1} ((0,T) \otimes R^d)
  $$
of two Banach spaces   $ L_{r_0, q_0} ((0,T) \otimes R^d)  $ and  $ L_{q_1, m_1} ((0,T) \otimes R^d)  $ is defined as a set
of all the functions  $ v = v(x,t), \ x \in R^d, \ t \in (0,T) $ of a form
$$
v(x,t) = v_0(x,t) + v_1(x.t),  \eqno(5.2)
$$
 where  $ v_0 \in L_{r_0, q_0}( (0,T) \otimes R^d), \ v_1 \in L_{q_1, m_1} ((0,T) \otimes R^d),  $ equipped with the norm

$$
|v|M_{r_0, q_0, q_1, m_1} = \inf \left[ |v_0| L_{r_0, q_0}( (0,T) \otimes R^d) + |v_1| L_{q_1, m_1} ((0,T) \otimes R^d)  \right], \eqno(5.3)
$$
 where $ "\inf" $ in (5.3) is calculated over all the functions $  v_0, v_1  $ satisfying the representation (5.2). \par
The space  $ M_{r_0, q_0, q_1, m_1}( (0,T) \otimes R^d) $ relative the norm (5.3) is also the complete Banach space. \par

It follows immediately from  theorems  2.0a and 3.1b the following assertion: \\

{\bf Proposition 5.2.}
 {\it  Let the tuple   }  $ r_0, q_0  $ {\it  satisfies the conditions (2.2) and
let the tuple  }  $ (q_1, m_1, k_1, h_1)  $ {\it  satisfies  the conditions (3.6). Then }

$$
|u|M_{r_0, q_0, q_1, m_1} \le  C_0(d; p_0) \ D(d,r_0) \ V(p_0,r_0) \ |f|_{p_0}   +
$$

$$
 D(d,r_1) \ K_B(1; r_1,k_1) \ C_0(h_1,r_1) \ V(h_1,r_1) \ |F|_{k_1,T;h_1,X}. \eqno(5.4)
$$

\vspace{4mm}

\section{Generalized Grand Lebesgue Spaces  estimations.}

\vspace{4mm}

 We assume in this section that the initial condition $ f = f(x) $ belongs to some Grand Lebesgue Space
$  G(\psi) =  G(\psi, X). $  This imply by definition that the following norm is finite:

$$
||f||G(\psi) := \sup_{ p \in (a,b) } \left[ \frac{|f|_p}{\psi(p)} \right] < \infty, \ 1 \le a < b \le \infty. \eqno(6.0)
$$
 Here $ \psi = \psi(p) $ is some positive continuous in {\it  open } interval $ (a,b) $ such  that \\
  $ \inf_{p \in (a,b)} \psi(p) > 0. $\par
 The detail investigation of these spaces see in \cite{Kozachenko1}, \cite{Fiorenza3}, \cite{Fiorenza4},
 \cite{Iwaniec1}, \cite{Iwaniec2}, \cite{Jawerth1}, \cite{Ostrovsky1}, \cite{Ostrovsky2} etc. \par

 Denote by $   M_0(p) = M_0(p; q,r), \ p \in (a,b), \ (q,r) \in G_0(p) $ the minimal value of the coefficient
in the inequality (2.3):

$$
  M_0(p; q,r):= \sup_{0 \ne f \in L_p(X)} \left[ \frac{ |u_0|_{r,T; q,X}}{ C_0(d; p) \ D(d,r) \ V(p,r) \ |f|_p} \right],  \eqno(6.1)
$$
we know that $ M_0(p; q,r) \le 1. $ Denote

$$
H_0 = \{  (q,r): \ \exists p \in (a.b) \ \Rightarrow (q,r) \in G_0(p); \}
$$

$$
\nu_0 =  \nu_0(d; r,q) = \inf_{p \in (a,b)} \left[ M_0(p; q,r) \ C_0(d; p) \ D(d,r) \ V(p,r) \ \psi(p)  \right], \ (q,r) \in H_0,
$$
and define the following  "two - dimensional"  generalization of the Grand Lebesgue Norm:

$$
||u_0||G \nu_0  \stackrel{def}{=}  \sup_{ (q,r) \in H_0} \left[ \frac{|u_0|_{r,T; q,X}}{ \nu( q,r) } \right]. \eqno(6.2)
$$

\vspace{3mm}

{\bf Theorem 6.1.}

$$
||u_0||G \nu_0 \le 1 \cdot ||f||G(\psi), \eqno(6.3)
$$
{\it  wherein the constant "1" in (6.3) is the best possible. }\par

\vspace{3mm}

{\bf Proof.} The {\it upper} estimate. Let $ f \in G(\psi);  $ we can and will conclude without loss of generality
$  ||f||G(\psi) = 1,  $ therefore

$$
|f|_p \le \psi(p), \ p \in (a,b).
$$
 It follows from theorem  2.0

$$
  |u_0|_{r,T; q,X} \le M_0(p; q,r) \ \psi(p),
$$
which is equivalent to (6.3) after minimisation over $  p. $  \par
 {\it The exactness } of the constant "1" is in fact proved in the article \cite{Ostrovsky208}, where is considered
the one - dimensional case; the multivariate version is completely analogous.\par

\vspace{3mm}

 Denote by $   M_1 = M_1(p, q,r,k), \ (q,r) \in G_1(p,q) $ the minimal value of the coefficient
in the inequality (3.2):

$$
  M_1(p, q,r,k):= \sup_{0 \ne F \in L_{r,X;k,T}} \left[ \frac{ |u_1|_{r,T; q,X}}{ K_1(p,q,r,k) \ |F|_{r,k}} \right],  \eqno(6.4)
$$
we know that $ M_1(p, q,r,k) \le 1. $

 Let $  Y = \{  r,k  \}, \ r,k > 1 $   be some open non-empty domain in the positive quarter plane $ R^2 $
and let $  \zeta = \zeta(r,k) $ be continuous in the set $  Y  $ positive function such that

$$
\inf_{ (r,k) \in Y } \zeta(r,k) > 0.
$$

 Define as before the two-dimensional Grand Lebesgue Norm

 $$
 || F(\cdot, \cdot) ||G\zeta = \sup_{(r,k) \in Y } \ \left[ \frac{|F|_{r,X; k,T}}{\zeta(r,k)}   \right]. \eqno(6.5)
 $$

Denote

$$
H_1 = \{  (p,q): \ \exists (r,k):  \ \Rightarrow (p,q,r,k) \in G_1(p,q) \};
$$

$$
\tau_1 =  \tau_1(d; p,q) = \inf_{(r,k) \in G_1(p,q)} \left[ M_1(p, q,r,k) \ \zeta(r,k)  \right], \ (q,r) \in H_1.
$$

\vspace{3mm}

{\bf Theorem 6.2.}

$$
||u_1||G \tau_1 \le 1 \cdot ||F||G(\zeta), \eqno(6.7)
$$
{\it  wherein the constant "1" in (6.6) is the best possible. }\par

\vspace{3mm}

{\bf Proof }  is at the same  as in theorem 6.1  through  theorem 3.1 and may be omitted.\\

\vspace{4mm}

\section{ Weight estimates. }

\vspace{4mm}

 Denote

 $$
 f(y) = g(y) \ |y|^{-a}, \ v_0 = v_0(x,t) = u_0(x,t) \ |x|^{-b} \ t ^{-\theta}, \ a, \ b  = \const \ge 0. \eqno(7.0)
 $$

 We try to obtain in this subsection first of all the {\it  weight }  estimations of a form

 $$
 | v_0|_{r,T; q,X} \le K_{a, b, \theta}(p,q,r) \ |g|_{p,X}, \eqno(7.1)
 $$
i.e. the weight generalization of theorem 2.0. \par

 {\it New notations and restrictions:}

$$
p_- = \frac{d}{d - a} \ge 1, \hspace{5mm} p_+ = \frac{d}{2/r - a - 2 \theta} > 0,
$$

$$
\kappa: = \frac{a + b + d - 2 \theta - 2/r}{d} \in (0,1),  \eqno(7.2)
$$

$$
p,q,r > 1,  \hspace{6mm}  p \in (p_-, \ p_+),
$$

$$
a + 2 \theta  < \frac{2}{r} < d + 2 \theta,
$$

$$
1+ \frac{1}{q}  = \frac{1}{p} + \kappa. \eqno(7.3)
$$

It is easy to see that the new designations agreed with the old in the case when $ a = b = 0. $ \par

\vspace{4mm}

{\bf Theorem 7.0.} {\it  We conclude under our restrictions (7.2) - (7.3) that the "constant" } $ K_{a, b, \theta}(p,q,r) $
{\it in  is finite:}

$$
  K_{a, b,\theta}(p,q,r)  \le \frac{C_2(d; a,b, \theta) \ D(d+ 2 \theta,r)}{ \left[p - p_- \right]^{\kappa} }. \eqno(7.4)
$$

 {\it Moreover, if }  $ K_{a, b, \theta}(p,q,r)  $ {\it is finite, then the triple $ (p,q,r) $ satisfies the equality (7.3)
 and $  p \in (p_-, \ p_+). $  }

\vspace{3mm}

{\bf Proof } is the same as one in theorem 2.0. We start from the inequality (2.4):

$$
 |v_0|_{r,T} \le |x|^{-b} \  \int_X |t^{-\theta} w_t(x-y)|_{r,T} \ |g(y)| \ |y|^{-a} \ dy =
$$

$$
 |x|^{-b} \ D(d+ 2 \theta,r) \int_X \frac{|g(y)| \ |y|^{-a} \ dy}{|x-y|^{d + 2\theta - 2/r} }=
D(d+ 2 \theta,r) \ I_{a, b, 2/r - 2 \theta}[g](x).  \eqno(7.5)
$$

 It remains to use again the main results of the articles \cite{Ostrovsky108}, \cite{Ostrovsky109}. \par

The {\it necessity}  of the conditions (7.2) - (7.3) may be proved as before  by means of the
scaling method.  By our opinion, the scaling method we used is somewhat simpler than in \cite{Kerman1}.  \par

 Note that the weight generalization of Young's convolution inequality  was considered at first by  R.A.Kerman in
\cite{Kerman1}; see also  the articles \cite{Nursultanov1},  \cite{Stepanov1},  \cite{Vasil'eva1} etc.

\vspace{4mm}

\section{Fractional Laplace operator.   }

\vspace{4mm}

 Let us consider the linear heat type initial value problem  with fractional power of Laplace operator:

 $$
 u_t + (-\Delta)^{\alpha}u = F(x,t), \ u(x,0+) = f(x), \ \alpha = \const >0, \  \alpha \ne 1,2,3. \ldots \eqno(8.0)
 $$
 We intend to generalize {\it weight} obtained estimations and some estimates of the works \cite{Miao1},  \cite{Miao2}, \cite{Zhai1}
 on the function $ u^{(\alpha)}(\cdot, \cdot) = u(\cdot, \cdot). $ \par

 Denote

 $$
 Z(x) = (2 \pi)^{-d/2} \int_{R^d} e^{i (x, \xi) - |\xi|^{2 \alpha} } d \xi, \ x \in R^d; \eqno(8.1)
 $$

$$
Z(x) = Z^{(0)}(|x|), \hspace{5mm} Z_t(x) = t^{-d/(2 \alpha)} \cdot Z \left( \frac{x}{t^{1/(2 \alpha)}}   \right), \ t > 0. \eqno(8.2)
$$

 The ordinary solution of equation (8.0)  has a view  as before

$$
 u^{(\alpha)}(x,t) = \int_X Z_t(x-y) \ f(y) \ dy + \int_0^{t} \int_X Z_{t-s}(x-y) \ F(y,s) \ ds  =
$$

$$
 u^{(\alpha)}_0(x,t) +   u^{(\alpha)}_1(x,t),  \hspace{5mm}  {\bf where} \hspace{5mm}   u^{(\alpha)}_0(x,t) = \int_X Z_t(x-y) \ f(y) \ dy, \eqno(8.3a)
$$

$$
u^{(\alpha)}_1(x,t) = \int_0^{t} \int_X Z_{t-s}(x-y) \ F(y,s) \ ds.  \eqno(8.3b)
$$

 Recall that the function

$$
q_d(x) = (2 \pi)^{-d} \int_X e^{i (x, \xi) - |\xi|^{2 \alpha} } d \xi
$$
is named {\it transstable density}, and appears in the Probability Theory,
see, e.g. \cite{Uchaikin1}, \cite{Fang1}. Evidently, $ Z(x) = (2 \pi)^{d/2} q_d(x, 2 \alpha). $ \par

 The asymptotic expression as $ |x| \to \infty $ for the transstable density is calculated in the book of Uchaikin V.V. and  Zolotarev V.M.
   \cite{Uchaikin1}, p. 212 - 215; see also \cite{Fang1}: as $ |x| \to \infty $

$$
|Z(x)| \sim C_{U, Z}(d,\alpha) \ (1 + |x|)^{-d - 2 \alpha}, \ x \in R^d. \eqno(8.4a)
$$

  The non-asymptotic  estimation for $ Z(x)  $ is calculated  by  Changxing Miao,  Baoquan Yuan,  and Bo Zhangin in \cite{Miao1}:

$$
|Z(x)| \le C_{M, Y, Z}(d,\alpha) \ (1 + |x|)^{-d - 2 \alpha}, \ x \in R^d. \eqno(8.4b)
$$

  Obviously, the function $  Z = Z(x) $ is radial, symbolically: $  Z = Z(|x|). $\par

  It follows from (8.4a) that the last estimation (8.4b) is asymptotically exact as $ |x| \to \infty $  up to multiplicative constant. \par

 We find by direct calculation from (8.4b) \hspace{5mm} $ |Z_t|_{p,X} \le  $

  $$
   C_{M, Y, Z}^{1/p}(d,\alpha) \ \omega^{1/p}(d) \ t^{d(1/p - 1)/(2 \alpha)} \ B^{1/p}(d - \alpha p, (p-1)d + 2 \alpha p) \le
 $$

$$
C_{M, Y, Z}^{1/p}(d,\alpha) \ \omega^{1/p}(d) \ t^{d(1/p - 1)/(2 \alpha)} \ B^{1/p}(d - \alpha, 2 \alpha), \ p \ge 1; \eqno(8.5)
$$

$$
d - \alpha p > 0, \hspace{5mm} (p-1)d + 2 \alpha p > 0;
$$

$$
|Z_t|_{q,T} \le (2 \alpha C_{M, Y, Z}(d,\alpha))^{1/q} \ |x|^{2 \alpha/q - d} \ B^{1/q}(d q - 2 \alpha; q(d + 2 \alpha) - d),\ q \ge 1,  \eqno(8.6)
$$
where $ B(\cdot, \cdot) $ is ordinary Beta - function and

$$
d q - 2 \alpha > 0,  q(d + 2 \alpha) - d) > 0,
$$
where
$$
\omega(d) = \frac{2 \ \pi^{d/2}}{\Gamma(d/2)}
$$
is an area of unit sphere in the Euclidean space $  R^d. $ \par

\vspace{4mm}

 We used the elementary integral

 $$
 \int_0^{\infty} \frac{x^{\gamma_1} \ dx }{(1+x)^{\gamma_2}} = B(\gamma_1 + 1, \gamma_2 - \gamma_1 - 1),   \eqno(8.7)
 $$

$$
\gamma_1, \gamma_2 = \const, \ \gamma_1 > -1, \ \gamma_2 - \gamma_1 > 1.
$$

 Moreover, let us denote also

$$
Z^{(\theta, \tau)}_t(x) = |x|^{-\theta} \ t^{-\tau} \ Z_t(x), \  \theta,\tau = \const.  \eqno(8.8)
$$

 We have using (8.7) the following weight estimates:

$$
|Z_t^{(\theta,\tau)}(\cdot)|_{r,T} \le (2 \alpha C_{M, Y, Z}(d + 2 \alpha \tau),\alpha))^{1/r} \ |x|^{2 \alpha/r - d - \theta} \times
$$

$$
B^{1/q}( (d + 2 \alpha \tau) q - 2 \alpha, q(d + 2 \alpha \tau + 2 \alpha) - d),\ r \ge 1,  \eqno(8.9)
$$
if of course

$$
(d + 2 \alpha \tau) r - 2 \alpha > 0, \ r(d + 2 \alpha \tau + 2 \alpha) - d > 0.  \eqno(8.10)
$$

\vspace{3mm}
 Analogously

 $$
 |Z_t^{\theta,\tau} (\cdot) |_{p,X} \le  C_{M, Y, Z}^{1/p}(d,\alpha) \ \omega^{1/p}(d) \ t^{-\tau + d(1/p - 1)/(2 \alpha) - \theta/(2 \alpha) } \times
 $$

 $$
   B^{1/p}(d - \alpha p, dp + 2 \alpha p - d + \theta p),  \eqno(8.11)
 $$
under restrictions

$$
d - \alpha p > 0, \ dp + 2 \alpha p - d + \theta p  > 0. \eqno(8.12)
$$

\vspace{3mm}

 Denote for brevity

$$
\lambda_1 =\lambda_1(r; d,\alpha, \theta, \tau) = - 2 \alpha/r + d + \theta,  \eqno(8.13)
$$

$$
\lambda_1^0 =\lambda_1(r; d,\alpha, \theta, 0),  \eqno(8.13a)
$$

$$
\lambda_2=  \lambda_2(p; d, \alpha, \theta, \tau)   = \tau - d(1/p - 1)/(2 \alpha) + \theta/(2 \alpha),  \eqno(8.14)
$$

$$
\lambda_2^0 =  \lambda_2(p; d, \alpha, 0, \tau)   = \tau - d(1/p - 1)/(2 \alpha).  \eqno(8.14a)
$$

\vspace{3mm}

 Then the estimates (8.9) and (8.11) may be rewritten in more simple form under restrictions (8.10) and (8.12)

$$
|Z_t^{(\theta,\tau)}(\cdot)|_{q,T} \le C_q  \ |x|^{-\lambda_1}, \hspace{5mm}
 |Z_t^{(\theta,\tau )} (\cdot) |_{p,X} \le  C_p  \ t^{-\lambda_2}. \eqno(8.15)
$$

{\bf Remark 8.1.}  Evidently, the estimates (8.15) are exact up to multiplicative constant on the {\it whole}
real axis (semi-axis):

\vspace{3mm}

$$
|Z_t^{(\theta,\tau)}(\cdot)|_{q,T} \ge \tilde{C}_q  \ |x|^{-\lambda_1}, \hspace{5mm}
 |Z_t^{(\theta,\tau )} (\cdot) |_{p,X} \ge  \tilde{C}_p  \ t^{-\lambda_2}. \eqno(8.15a)
$$

\vspace{3mm}

{\bf Remark 8.2.}  More surprisingly is the circumstance that the estimates (8.15) and (8.15a) are true even
for the integer values $ \alpha, $ for instance for the value $ \alpha = 1, $ in which $ Z_t(x) = w_t(x). $\par

 This is all the more surprising that the tail behavior as $ |x| \to \infty $ of the function $ w_1(x)  $
 strikingly different on the ones for the function $ Z_1(x). $\par

\vspace{3mm}

 We note passing to the descriptions of next results first of all  that the
decay estimates as $  t \to \infty $ for solution $ u = u(x,t) $ of linear and nonlinear nonlocal heat equations, for more general
pseudo-differential  equation, of a form

$$
|u|_{p,X} \le C \ t^{-[(d/\gamma) \cdot (1/q - 1/p)] } \ |f|_{q,X},  \ 1 \le q < p, \ \gamma = \const \in (0,\infty),
$$
(in our notations),   was obtained in recent  article  \cite{Brandle1}; see also reference therein. \par

\vspace{3mm}

 Let us consider now the inverse order of the spaces.  A new notations and restrictions:

 $$
 v^{(\alpha)}_0(x,t) =  v_0(x,t) = u^{(\alpha)}_0(x,t) \ |x|^{-b} \ t^{-\tau}, \ f(x) = g(x) \ |x|^{-a},
 $$

$$
a,b,\tau = \const \ge 0, \ \kappa_0 := (a + b + \lambda_1^0)/d  \in (0,1); \eqno(8.16)
$$

$$
p_- = \frac{d}{d-a} \ge 1, \ p_+ = \frac{d}{d - a - \lambda_1^0} >0. \eqno(8.17)
$$

$$
d + \frac{d}{q} = \frac{d}{p}+ (a + b + \lambda_1^0), \ q = q(p).  \eqno(8.18)
$$

\vspace{3mm}

{\bf Theorem 8.0.} {\it   We conclude under conditions (8.10), (8.16), (8.17) and (8.18) for the values $ p $ from the interval  }
$ p \in (p_-, p_+)  $  {\it  and for all the positive values $  \alpha,  $ including the integer numbers, for example $ \alpha = 1, $
in which $ Z_t(x) $ is the classical solution of the heat equation, }

$$
|v_0^{(\alpha)}|_{r,T; q,X} \le \frac{C(d; a,b,\tau,\alpha)}{[(p - p_-  ) \ (p_+ - p)]^{\kappa_0} } \ |f|_p,  \eqno(8.19)
$$
{\it  if, of course,  }  $ f \in L_p(R^d). $ \par
 {\it  Conversely, the equality (8.18) is necessary for the estimation of the form (8.19). }\par
{\it  Moreover, the estimation (8.19) is  exact up to multiplicative constant. } \par

\vspace{3mm}

{\bf Proof.}

$$
 v_0^{(\alpha)}(x,t) =  |x|^{-b} \ \int_X t^{-\tau} \ Z_t(x-y) \ g(y) \ |y|^{-a} \ dy =
$$

$$
|x|^{-b} \ \int_X  Z_t^{(0,\tau)}(x-y) \ g(y) \ |y|^{-a} \ dy.
$$

We have using again  Marcinkiewicz integral-triangle  inequality:

$$
|v_0^{(\alpha)}(x,t)|_{r,T} \le
|x|^{-b} \ \int_X | Z_t^{(0,\tau)}(x-y)|_{r,T} \ |g(y)| \ |y|^{-a} \ dy.
$$
 We use now the estimate (8.15):

$$
 |v_0^{(\alpha)}(x,t)|_{r,T} \le
|x|^{-b} \ \int_X C_p \ \frac{|y|^{-a} \ |g(y)| \ dy}{ |x-y|^{\lambda_1^0 }}. \eqno(8.20)
$$

 The integrals of a form (8.20) are called weight Hardy-Littlewood operators or equally weight fractional integral operators.
The using for us {\it exact } up to multiplicative constant  $  L_q(R^d) $ estimates for
them  are obtained in the article  \cite{Ostrovsky108}. \par

 The necessity of the equality (8.18) for the assertion of considered theorem   may  be proved as ordinary by means of scaling method.\par

This completes the proof of theorem (8.1). \par

\vspace{3mm}

 We will deduce now he weight estimates for the second term in the representations (8.3a), i.e. for the function $ u_1^{(\alpha)}(x,t). $ \par

\vspace{3mm}

 Some new notations and restrictions.

 The first  {\it weight} function $  v^{(\alpha)}_1(x,t) $ for the part of solution
 $  u^{(\alpha)}_1(x,t) $ is defined as follows:

$$
v^{(\alpha)}_1(x,t) := |x|^{-b} \ u^{(\alpha)}_1(x,t) = |x|^{-b} \int_0^t \int_X Z(x-y,t-s) \ |y|^{-a} \ H(y,s) \ ds \ dy, \eqno(8.21a)
$$
where  $ H(y,s) = |y|^{a} \cdot F(y,s) $ is the first weight function for the right - hand  side $ F = F(x,t). $  Here $ a,b = \const \ge 0. $\par

\vspace{3mm}

 The second  {\it weight} function $  v^{(\alpha)}_2(x,t) $ for the part of solution
 $  u^{(\alpha)}_1(x,t) $ is defined as follows:

$$
v^{(\alpha)}_2(x,t) := t^{-\tau} \ u^{(\alpha)}_1(x,t) = t^{-\tau} \int_0^t \int_X Z(x-y,t-s) \ s^{-\beta} \ R(y,s) \ ds \ dy, \eqno(8.21b)
$$
where  $ R(y,s) = s^{-\beta} \cdot F(y,s) $ is the second weight function for the right - hand  side $ F = F(x,t). $  Here $ \tau, \beta = \const \ge 0. $\par

 Other notations  and restrictions.

$$
1 + \frac{1}{r_1} =  \frac{1}{Q_1} + \frac{1}{m_1}, \hspace{6mm} r_1,Q_1,m_1 \in (1, \infty), \eqno(8.22a)
$$

$$
a,b = \const  \ge 0, \hspace{6mm} \kappa_1 = \frac{a + b + d - 2 \alpha}{d} \in (0,1),
$$

$$
1 + \frac{1}{q_1} =  \frac{1}{p_1} + \frac{a + b + d - 2 \alpha}{d} = \frac{1}{p_1} + \kappa_1 , \hspace{6mm} p_1, q_1 \in (1, \infty),
\eqno(8.22b)
$$

$$
q_1 = q_1(p_1), \   p_-^1 = \frac{d}{d - a} > 1, \ \hspace{4mm} p_+^1 = \frac{d}{2 \alpha/Q_1 - a} > 0,
$$

$$
p_1 \in (p_-^1, p_+^1).
$$

\vspace{3mm}

{\bf Theorem 8.1a} {\it  We propose under formulated assumptions  }

\vspace{3mm}

$$
|v^{(\alpha)}_1(\cdot,\cdot) |_{r_1,T;  q_1,X }   \le \frac{C_1(d;a,b)}{ [ (p_1 - p_-^1)(p_+^1 - p_1)  ]^{\kappa_1}} \ |H(\cdot, \cdot)|_{m_1,T; p_1,X}.
\eqno(8.23)
$$

{\it  Conversely, the estimation  (8.23)  is exact up to multiplicative constant and the relations (8.22a), (8.22b) are necessary
 for the inequality  of a form (8.23).  }\par

\vspace{3mm}

{\bf Proof.} We have from (8.21) using again the  Marcinkiewicz integral  triangle  inequality and  Beckner's convolution  constants:

$$
|v^{(\alpha)}_1(x,\cdot) |_{r_1,T}   \le |x|^{-b}  \int_X |y|^{-a} \ K_B(1; Q,m) \  |Z(x-y, t - s)|_{Q_1,T} \ |H|_{m_1,T} \ dy =
$$

$$
C \cdot K_B(1; Q_1,m_1) \cdot \int_X \frac{|x|^{-a} \ |y|^{-b} }{|x-y|^{d - 2 \alpha/Q_1}} \ |H(y,\cdot)|_{m_1,T} \ dy.
$$

 It remains to apply the inequality (2.4) for the weight Hardy - Riesz transform. \par
 The exactness of the inequality (8.23) follows immediately from the main result of the report \cite{Ostrovsky108}  by means of consideration
of the factorable function  $ F(x,t) = F_0(t) \  F_1(x).$ \par

  The necessity of the  conditions (8.22a) and (8.22b) for (8.23) may be  elementary provided  by the known scaling method, where the dilation
operator $ T_{\lambda, \mu}    $ is defined as follows

$$
T_{\lambda, \mu} [F_0 \  F_1](t,x) \stackrel{def}{=} F_0(\lambda t) \ F_1(\mu x),  \ \lambda, \mu \in (0, \infty).
$$

\vspace{4mm}

 Let us consider the inverse order of norms. Notations and restrictions: \par

$$
1 + \frac{1}{q_2} = \frac{1}{Q_2} +  \frac{1}{p_2},  \eqno(8.24a)
$$

$$
0 \le \beta, \tau < 1,  \hspace{4mm} \zeta_2:= d(1 - 1/Q_2)/(2 \alpha) \in(0,1),  \hspace{4mm} \beta + \tau + \zeta_2 < 1,
$$

$$
1 + \frac{1}{r_2} = \frac{1}{k_2} +  (\beta + \tau + \zeta_2), \ r_2 = r_2(k_2), k_2, r_2, q_2, Q_2, p_2 \in (1,\infty),
\eqno(8.24b)
$$

$$
\kappa_2 := \beta + \tau + \zeta_2, \ p_-^{(2)} = \frac{1}{1-\beta}, \ p_+^{(2)} = \frac{1}{1 - \beta - \tau} > 0,
$$

$$
p_2 \in( p_-^{(2)}, p_+^{(2)}).
$$

\vspace{4mm}

{\bf Theorem 8.1b} {\it  We propose under formulated assumptions  }

\vspace{3mm}

$$
|v^{(\alpha)}_2(\cdot,\cdot) |_{q_2,X; r_2,T }   \le \frac{C_2 K_B(d; Q_2, p_2)}{ (p_2 - p_-^{(2)} )^{\kappa_2} } \ |R(\cdot,s)|_{p_2,X; k_2,T }.
\eqno(8.25)
$$

{\it  As before, the estimation  (8.25)  is exact up to multiplicative constant and the relations (8.24a), (8.24b) are necessary
 for the inequality  of a form (8.25).  }\par

\vspace{3mm}

{\bf Proof.}  Note that

$$
v^{(\alpha)}_2(x,t) = t^{-\tau} \int_0^t s^{-\beta} \ ds \ \int_X Z(x-y, t - s) \ R(y,s) \ dy =
$$

$$
t^{-\tau} \int_0^t s^{-\beta} \ ds  \cdot Z_{t-s}*[R](x).
$$
 Therefore

 $$
| v^{(\alpha)}_2(x,t) |_{q_2,X} \le t^{-\tau} \int_0^t s^{-\beta} \ ds \cdot |Z_{t-s}*[R](x)|_{q_2,X} \le
 $$

$$
C_2 \ t^{-\tau} \int_0^t s^{-\beta} \ ds \cdot K_B(d; Q_2, p_2) \ (t-s)^{-d(1-1/p)/(2 \alpha)} \ |R(\cdot,s)|_{p_2,X}.
$$

 The proposition of theorem 8.1b follows immediately from the main result of an article
\cite{Ostrovsky109}. \par

\vspace{4mm}

\section{Concluding remarks.}

\vspace{3mm}

{\bf A. General  elliptic operator.   }\\

\vspace{3mm}

  Our results may be  easily generalized on the initial value problem for parabolic PDE of a form

$$
\frac{\partial u}{\partial t}  = 0.5 \sum_{k=1}^d \ \sum_{i=1}^d a_{i,k}(x,t) \frac{\partial^2 u }{ \partial x_i \partial x_k} +
\sum_{i=1}^d b_i(x,t) \frac{\partial u}{\partial x_i} + c (x,t) u + F(x,t), \eqno(9.1)
$$

$$
 u(x, 0+) = f(x), \hspace{6mm} x \in R^d, \ t > 0, \eqno(9.2)
$$
if for example the functions $ a_{i,k}(x,t)  $ are symmetrical $ a_{i,k}(x,t) = a_{k,i}(x,t),  $ uniform positive definite:

$$
\sum_{k=1}^d \sum_{i=1}^d a_{i,k}(x,t) \xi_i \xi_k \ge \lambda \sum_{k=1}^d \xi_k^2, \hspace{5mm} \lambda = \const > 0;
$$
all the functions $ a_{i,k}(x,t), \ b_i(x,t),  \ c(x,t)   $ are bounded and satisfy the H\"older's conditions with positive power. \par
 It is well-known that  under these conditions the solution of problem (9.1)-(9.2) may be written as follows:

 $$
 u(x,t) = \int_{R^d} \int_0^t G(x,t,y,s) \ F(y,s) \ dy \ ds + \int_{R^d} G(x,t,y,0) \ f(y) \ dy,
 $$
where the function $ G(x,t,y,s) $ allows an estimation

$$
G(x,t,y,s) \le C_1 (t-s)^{-d/2} \ \exp \left( - C_2 |x-y|^2/(t-s)  \right), \ 0 \le s < t < \infty.
$$

\vspace{3mm}

{\bf B. Arbitrary convolution kernel.   }\\

\vspace{3mm}

 Our considerations may be easily generalized on the function of a form

$$
z(x,t) = \int_X M(x-y,t) \ f(y) \ dy  + \int_0^t \int_X M(x-y, t-s) \ F(y,s) \ ds  \  dy,
$$
where the kernel function $  M(x,t) $ has the form

$$
M(x,t) = t^{-a} M_1(|x|/t^b), \  M_1(\cdot) \in  L_p(R^d).
$$

\end{document}